\documentclass[12pt,oneside,reqno]{amsart}
\usepackage{mathrsfs}
\usepackage{graphics}
\usepackage{graphicx,color}
\usepackage{amssymb}
\usepackage{subfigure}
\usepackage{enumerate}
\newcommand{\dif}{\mathrm{d}}

\newcommand{\be}{\begin{eqnarray}}
\newcommand{\ee}{\end{eqnarray}}
\newcommand{\ce}{\begin{eqnarray*}}
\newcommand{\de}{\end{eqnarray*}}
\newtheorem{theorem}{Theorem}[section]
\newtheorem{lemma}[theorem]{Lemma}
\newtheorem{remark}[theorem]{Remark}
\newtheorem{definition}[theorem]{Definition}
\newtheorem{proposition}[theorem]{Proposition}
\newtheorem{Example}[theorem]{Example}
\newtheorem{corollary}[theorem]{Corollary}
\def\e{\varepsilon}
\def\t{\theta}
\def\a{\alpha}

\def\[{{\Big[}}
\def\]{{\Big]}}
\def\<{{\langle}}
\def\>{{\rangle}}
\def\({{\Big(}}
\def\){{\Big)}}

\def\no{\nonumber}
\def\bt{\begin{theorem}}
\def\et{\end{theorem}}
\def\bl{\begin{lemma}}
\def\el{\end{lemma}}
\def\br{\begin{remark}}
\def\er{\end{remark}}
\def\bx{\begin{Example}}
\def\ex{\end{Example}}
\def\bd{\begin{definition}}
\def\ed{\end{definition}}
\def\bp{\begin{proposition}}
\def\ep{\end{proposition}}
\def\bc{\begin{corollary}}
\def\ec{\end{corollary}}

\def\cB{{\mathcal B}}
\def\cC{{\mathcal C}}

\def\cI{{\mathcal I}}
\def\cJ{{\mathcal J}}

\def\cM{{\mathcal M}}
\def\cN{{\mathcal N}}

\def\mE{{\mathbb E}}

\def\mQ{{\mathbb Q}}
\def\mR{{\mathbb R}}

\def\mU{{\mathbb U}}

\def\geq{\geqslant}
\def\leq{\leqslant}

\def\sF{{\mathscr F}}

\def\sU{{\mathscr U}}

\title{Effective Filtering for Multiscale Stochastic Dynamical Systems driven by L\'evy processes*}

\author{Huijie Qiao}

\dedicatory{School of Mathematics, Southeast University\\
Nanjing, Jiangsu 211189,  China\\
hjqiaogean@seu.edu.cn}

\thanks{{\it AMS Subject Classification(2010):} 60H10; 37D10, 70K70, 60G51.}

\thanks{{\it Keywords:} Multiscale systems driven by L\'evy processes; random slow manifolds; dimension reduction; efficient filtering.}

\thanks{*This work was supported by NSF of China (No. 11001051, 11371352) and China Scholarship Council under Grant No. 201906095034.}

\subjclass{}

\date{}

\begin{document}

\allowdisplaybreaks

\begin{abstract}
The work is about multiscale stochastic dynamical systems driven by L\'evy processes. First, we prove that these systems can approximate low-dimensional systems on random invariant manifolds. Second, we establish that nonlinear filterings of multiscale stochastic dynamical systems also approximate that of reduced low-dimensional systems. Finally, we investigate the reduction for $\e=0$ and obtain that these reduced systems does not approximate these multiscale stochastic dynamical systems.
\end{abstract}

\maketitle \rm

\section{Introduction}

Given a probability space $(\Omega, \sF, \mQ)$ (See Subsection \ref{mds} in detail). Consider the following stochastic slow-fast system on $\mR^n\times\mR^m$:
\be\left\{\begin{array}{l}
\dot{u}^\e=\frac{1}{\e}A{u}^\e+\frac{1}{\e}U(u^\e,v^\e)+\frac{\sigma_1}{\e^{1/\alpha}}\dot{L}^{\alpha\pm}, \\
\dot{v}^\e=B{v}^\e+V(u^\e, v^\e)+\sigma_2\dot{L}^\pm,
\label{slfasy}
\end{array}
\right.
\ee
where $A$ and $B$ are $n\times n$ and $m\times m$ matrices respectively,  and the interaction
functions $U:\mR^n\times\mR^m\mapsto \mR^n$ and $V:\mR^n\times\mR^m\mapsto \mR^m$ are
Borel measurable. $L^{\alpha\pm}$ and $L^\pm$ are a $n$-dimensional two-sided symmetric $\a$-stable process and a $m$-dimensional two-sided L\'evy process, respectively(c.f. Subsection \ref{lepro}). Moreover, $\sigma_1$ and $\sigma_2$ are nonzero real noise intensities, and $\e$ is a small positive parameter representing the ratio of the two time scales. These systems like (\ref{slfasy}) are usually called multi-scale systems and have been applied to simulate many phenomena in chemistry, biology, climate, and so on(\cite{mg, Tur, wu, KP, ZS, zhang}). 

If L\'evy processes are replaced by Brownian motions, these systems (\ref{slfasy}) have been widely observed and studied. Let us recall some referrences. Khasminskii and Yin \cite{RG} developed a stochastic averaging principle for systems (\ref{slfasy}). Later, Schmalfu{\ss} and Schneider \cite{Sch} and Wang and Roberts \cite{WR} obtained the invariant manifold for systems (\ref{slfasy}).  In the infinite dimensional framework Fu, Liu and Duan \cite{Fu} and the author \cite{q3} studied the invariant manifolds of systems (\ref{slfasy}) and obtained low dimensional reduced systems. Besides, if systems (\ref{slfasy}) are driven by two symmetric $\a$-stable processes, i.e. the L\'evy process $L$ is also a symmetric $\a$-stable process, the author and two coauthors \cite{zqd} showed the existence of  the invariant manifold for systems (\ref{slfasy}). 

The nonlinear filtering  problem for systems (\ref{slfasy}) with respect to an observation process $\{w^{\e}_s, 0\leq s\leq t\}$ (See Subsection \ref{nonfil} in details)  is to evaluate the `filter' $\mE[\Phi(u^{\e}_t, v^{\e}_t)|\mathcal{W}^{\e}_t]$, where $\Phi$ is a Borel measurable function such that $\mE|\Phi(u^{\e}_t, v^{\e}_t)|<\infty$ for $t\in[0,T]$, and $\mathcal{W}^{\e}_t$ is the $\sigma$-algebra generated by $\{w^{\e}_s, 0\leq s\leq t\}$ (\cite{BC, RBL}). And it is sometimes called data assimilation (\cite{mg, ge, zhang}). Moreover, the nonlinear filtering problems for a number of multi-scale systems have been alternatively investigated (\cite{Imkeller} \cite{mg}-\cite{Park3} \cite{q2}-\cite{q4} \cite{zhang} \cite{zqd}). Let us mention some results. With the help of stochastic averaging, Park and his coauthors \cite{Park1}-\cite{Park3} studied filtering problems for two time scales systems with Brownian motions by the Zakai equations. And then Imkeller et al. \cite{Imkeller} showed that the filtering of only the slow part converges to the homogenized filtering by double backward stochastic differential equations and asymptotic techniques when the system (\ref{slfasy}) is driven by Brownian motions. Recently, the author and two coauthors \cite{q4} reduced the system (\ref{slfasy}) with Brownian motions to a system on a random invariant manifold, and showed that  the filtering of only the slow part converges to the filtering of the reduced system. Later, the author and two coauthors \cite{zqd} extended the result in \cite{q4} to the case of two symmetric $\a$-stable processes.

In the paper, we consider the system (\ref{slfasy}), where $L^\pm$ can be not only a two-sided symmetric $\a$-stable process but also a general L\'evy process. First, it is proved that these systems can approximate low-dimensional systems on random invariant manifolds. Second, we establish that nonlinear filterings of multiscale stochastic dynamical systems, rather than only the slow part, also approximate that of those reduced low-dimensional systems. Third, we deduce the reduced system for $\e=0$ and find that the system (\ref{slfasy}) does not converge to the reduced system as $\e\rightarrow0$.

Here our motivation is three-folded. The first fold is to correct the estimate for the distance between the system (\ref{slfasy}) and the reduced system in \cite[Theorem 3.2]{q4} and \cite[Theorem 2]{zqd}. That is, the estimate should {\it not} just depend on the initial value of the system (\ref{slfasy}) since the reduced system is on an invariant manifold (See Theorem \ref{resyth} in detail). The second fold is to extend the result in \cite[Theorem 3]{zqd}. Since symmetric $\a$-stable processes are a special type of L\'evy processes, the extension to general L\'evy processes is necessary. Finally, we analysis the case for $\e=0$. It is unfortunate to obtain that the system (\ref{slfasy}) does not approximate the reduced system as $\e\rightarrow0$.

It is worthwhile to mention our methods. First of all, note that our conditions are similar to those in \cite{Sch}. There Schmalfu{\ss} and Schneider made complicated deduction in order to construct an invariant manifold and furthermore implicitly expressed the manifold. Here we construct an invariant manifold only by an integral equation. Besides, to obtain the reduced system on the invariant manifold, we only define an operator and then prove that it is contractive. Thus, a large number of computation like that in \cite{Fu} is avoided.

\bigskip

This paper is arranged as follows. In Section \ref{pre},  we introduce basic concepts of random dynamical systems, stationary solutions, random invariant manifolds  and L\'evy processes. The existence of low-dimensional systems approximating these multiscale systems is placed in Section \ref{resyinma}. In Section \ref{filter}, we  introduce nonlinear filtering problems and prove that the nonlinear filterings of the low-dimensional reduced systems also approximate that of the multiscale systems. Next, we analysis the case for $\e=0$ in Section \ref{epze}. In Section \ref{con}, we summarize all the results in the paper. 

The following convention will be used throughout the paper: $C$ with or without indices will denote
different positive constants whose values may change from one place to
another. 

\section{Preliminaries}\label{pre}

In the section, we introduce basic concepts of random dynamical systems, stationary solutions, random invariant manifolds  and L\'evy processes.

\subsection{Random dynamical systems (\cite{la})}

\bd
Let $(\Omega,\sF,\mQ)$ be a probability space, and $(\theta_t)_{t\in\mR}$ a family of measurable
transformations from $\Omega$ to $\Omega$.
We call $(\Omega,\sF,\mQ; (\theta_t)_{t\in\mR})$ a metric dynamical system if for each $t\in\mR$,
$\theta_t$ preserves the probability measure $\mQ$, i.e.,
$$
\theta_t^*\mQ=\mQ,
$$
and for $s,t\in\mR$,
$$
\theta_0=1_\Omega, \quad \theta_{t+s}=\theta_t\circ\theta_s.
$$
\ed
\bd\label{rds} Let $(\mU,\sU)$ be a measurable space. A mapping
\ce
\Psi: \mR\times\Omega\times\mU\mapsto\mU, \quad
(t,\omega,x)\mapsto\Psi(t,\omega,x)
\de
with the following properties is called a measurable random dynamical system (RDS in short), or a cocycle:

(i) Measurability: $\Psi$ is
$\cB(\mR)\otimes\sF\otimes\sU/\sU$-measurable,

(ii) C\`adl\`ag cocycle: $\Psi(t,\omega)$ is c\`adl\`ag
for $t\in\mR$,  and further satisfies the following conditions
\be
\Psi(0,\omega)&=&id_{\mU},  \label{perfect coc1}\\
\Psi(t+s,\omega)&=&\Psi(t,\theta_s\omega)\circ\Psi(s,\omega), \label{perfect coc2}
\ee
for all $s,t\in\mR$ and $\omega\in\Omega$.
\ed

\subsection{\textbf{Stationary solutions and random invariant manifolds} (see \cite{la})} Let $\Psi$ be a RDS on the normed space $(\mU, \|\cdot\|_{\mU}, \sU)$, where $\|\cdot\|_{\mU}$ stands for a norm on $\mU$ and $\sU$ is Borel $\sigma$-field on $\mU$. We introduce stationary solutions and random invariant manifolds with respect to $\Psi$.

A random variable $\zeta$ is called a {\it stationary solution} of a stochastic/random differential
equation, if $\Psi$ is defined by the solution mapping of the equation and for $t>0$
$$
\Psi(t,\omega,\zeta(\omega))=\zeta(\theta_t\omega),\  a.s.\omega.
$$
Here, we remind that in general it is not obvious whether the solutions to stochastic differential equations define RDSs. Therefore, before mentioning stationary solutions, we need to justify that solution mappings of stochastic differential equations define RDSs.

A family of nonempty closed sets $\mathcal{M}=\{\mathcal{M}(\omega)\}_{\omega \in \Omega}\subset\sU$ is called a random set if for every $ u \in\mU $, the mapping
$$
\Omega\ni \omega\mapsto{\rm dist}(u,\mathcal{M}(\omega)):=\inf\limits_{u^{'}\in \mathcal{M}(\omega)}||u-u^{'}||_{\mU}
$$
is measurable. Moreover, $ \mathcal{M} $ is called a positively invariant set with respect to the random dynamical system $\Psi$ if
\begin{equation}\label{invariant set}
\Psi(t,\omega,\mathcal{M}(\omega))\subseteq \mathcal{M}(\theta_t\omega), \, \text{ for }\ t\in\mathbb{R^{+}},\ \omega\in\Omega.
\end{equation}

In the sequel, we consider a random set defined by a Lipschitz continuous graph. Concretely speaking, we define a function by $$
\Omega\times\mR^m \ni (\omega,y)\mapsto F(\omega, y)\in\mR^n
$$
such that for all $\omega\in\Omega$, $F(\omega, y)$ is globally Lipschitzian in $y$  and for any
$y\in\mR^m$, the mapping $\omega\rightarrow F(\omega, y)$ is a random vector.  And set
$$
\cM(\omega) :=\{(F(\omega, y), y), y\in\mR^m\},
$$
 and then $\cM$ is a random set (\cite[Lemma 2.1]{Sch}). Moreover, the random set
$\cM(\omega)$ is called a {\it Lipschitz random invariant manifold} if it is (positively) invariant with respect to some random dynamical system.

% or {\it Lipschitz random slow manifold}.

\subsection{\textbf{L\'evy processes }(see \cite{sa})}\label{lepro}
\bd\label{levy1} A stochastic process $L=(L_t)_{t\geq0}$ with $L_0=0$ a.s. is a
$n$-dimensional L\'evy process if

(i) $L$ has independent increments; that is, $L_t-L_s$ is
independent of $L_v-L_u$ if $(u,v)\cap(s,t)=\emptyset$;

(ii) $L$ has stationary increments; that is, $L_t-L_s$ has the same
distribution as $L_v-L_u$ if $t-s=v-u>0$;

(iii) $L_t$ is right continuous with left limit. \ed

Its characteristic function is given by \ce
\mE\big(\exp\{i\<z,L_t\>\}\big)=\exp\{t\varphi(z)\}, \quad z\in\mR^n.
\de The function $\varphi: \mR^n\rightarrow\mathcal {C}$ is called the
characteristic exponent of the L\'evy process $L$. By the
L\'evy-Khintchine formula, there exist a nonnegative-definite
$n\times n$ matrix $Q$, $b\in\mR^n$ and a measure $\nu$ on $\mR^n$
satisfying
\be
\nu(\{0\})=0 ~\mbox{and}~ \int_{\mR^n\setminus\{0\}}(|u|^2\wedge1)\nu(\dif u)<\infty, \label{lemc}
\ee
such that
\be
\varphi(z)=-\frac{1}{2}\<z,Qz\>+i\<z,b\>
+\int_{\mR^n\setminus\{0\}}\big(e^{i\<z,u\>}-1-i\<z,u\>1_{|u|\leq\delta}\big)\nu(\dif
u), \label{lkf} 
\ee 
where $\delta$ is a positive constant. $\nu$ is called the L\'evy measure associated with $L$.

Set $\kappa_t:=L_t-L_{t-}$. Then $\kappa$ defines a stationary Poisson point process with values in
$\mR^n\setminus\{0\}$ and the characteristic measure $\nu$ (\cite{iw}). Let $N_{\kappa}((0,t],\dif u)$ be the counting measure
of $\kappa_{t}$, i.e., for $D\in\cB(\mR^n\setminus\{0\})$
$$
N_{\kappa}((0,t],D):=\#\{0<s\leq t: \kappa_s\in D\},
$$
where $\#$ denotes the cardinality of a set. The compensated martingale measure
of $N_{\kappa}$ is given by
$$
\tilde{N}_{\kappa}((0,t],\dif u):=N_{\kappa}((0,t],\dif u)-t\nu(\dif
u).
$$
The L\'evy-It\^o theorem states that there exist a $n'$-dimensional Brownian motion $R_t$, $0\leq n'\leq n$ and  a
$n\times n'$ matrix $M$ such that $L$ can be represented as
\ce
L_t=bt+MR_t+\int_0^t\int_{|u|\leq\delta}u\tilde{N}_{\kappa}(\dif s, \dif u)
+\int_0^t\int_{|u|>\delta}uN_{\kappa}(\dif s, \dif u).
\de

In the sequel, we take a $m$-dimensional L\'evy process
\ce
L_t^+=MR_t+\int^t_0\int_{|u|\leq\delta}u\,\tilde{N}_\kappa(\dif s,\dif u), ~~\qquad\qquad t\geq0.
\de
Here, for convenience of the following deduction, we require $0<\delta<1$. And then set
\be\label{lif1}
L_t^\pm:=\left\{\begin{array}{l}
L_t^+, ~~~\qquad\quad t\geq0,\\
-L_{(-t)-}^+, ~\quad t<0.
\end{array}
\right.
\ee
Thus, $L^\pm$ is a two-sided $m$-dimensional L\'evy process. 

\vspace{3mm}

\begin{definition}
For $\alpha \in (0,2)$, a $n$-dimensional symmetric $\alpha $-stable process $ L^{\alpha}_{t} $ for $t\geq 0$ is a L\'evy process with the characteristic exponent 
\ce
\varphi(u)=-C_1(n,\alpha)| u |^{\alpha},  ~for~u \in {\mathbb{R}^{n}},
\de
where $ C_1(n, \alpha):=\pi^{-\frac{1}{2}}\Gamma((1+\alpha)/2)\Gamma(n/2)/\Gamma((n+\alpha)/2)$.
\end{definition}

For a $n$-dimensional symmetric $\alpha$-stable L\'evy process, the diffusion matrix $ Q= 0$,
the drift vector $ b= 0$, and the L\'evy measure $\nu $ is given by
\ce
\nu(du)=\frac{C_2(n,\alpha)}{{| u |}^{n+\alpha}}du,
\de
where $ C_2(n, \alpha):=\alpha\Gamma((n+\alpha)/2)/{(2^{1-\alpha}\pi^{n/2}\Gamma(1-\alpha/2))}$.

In the sequel, we fix a $n$-dimensional symmetric $\alpha $-stable process $ L^{\alpha}_{t} $ ($1<\alpha<2$) for $t\geq 0$ independent of $L^\pm=(L_t^\pm)_{t\in\mR}$ and then set
\be\label{alpha}
L^{\alpha\pm}_t:=\left\{\begin{array}{l}
L_t^{\alpha}, ~~~\qquad\quad t\geq0,\\
-L_{(-t)-}^{\alpha}, ~\quad t<0.
\end{array}
\right.
\ee
Thus, $L^{\alpha\pm}_t$ is a two-sided $n$-dimensional symmetric $\alpha$-stable process.

\section{The reduction system on a random invariant manifold}\label{resyinma}

In the section, we prove that a fast-slow system driven by L\'evy processes can approximate a low dimensional system on a random invariant manifold.

\subsection{A metric dynamical system}\label{mds}

Let $D(\mR,\mR^n)$ be the set of all functions which are c\`adl\`ag for $t\in\mR$, and take values in $\mR^n$. We take the canonical sample space $\Omega^1 \triangleq  D(\mR,\mR^n)$. It,  endowed with the Skorohod metric $\rho$, can be made a complete and separable metric space. The Borel $\sigma$-algebra on the sample space $\Omega^1$  under the
topology induced by $\rho$ is denoted as $\sF^1$. Let $\mQ^1$ be the distribution of the two-sided $n$-dimensional symmetric $\alpha$-stable L\'evy process $L^{\alpha\pm}=(L^{\alpha\pm}_t)_{t\in\mR}$. Set
\ce
\begin{aligned}
  &\theta^1:\mathbb{R}\times\Omega^1\mapsto \Omega^1,\\
  &\theta^1_{t}\omega(\cdot):=\omega(\cdot+t)-\omega(t),
  \end{aligned}
\de
and then one can justify that the probability measure $ \mathbb{Q}^1$ is $ \theta^1$-invariant and $\{\theta^1_{t}, t\in\mR\}$ is a group. Thus, $(\Omega^1,\sF^1,\mathbb{Q}^1, (\theta_t^1)_{t\in\mR})$ is a metric dynamical system. 

Next let $\Omega^2 \triangleq  D(\mR,\mR^m)$. Likewisely, we define $ \sF^2$ and $\theta_t^2 $. Again let $\mQ^2$ be the unique probability measure which makes the canonical process the L\'evy process $L^\pm=(L_t^\pm)_{t\in\mR}$. So, $ (\Omega^2, \sF^2,\mathbb{Q}^2, (\theta_t^2)_{t\in\mR}) $ is another metric dynamical system.

 Set
\ce
\Omega:=\Omega^1\times\Omega^2, ~\sF:=\sF^1\times\sF^2, ~\mathbb{Q}:=\mathbb{Q}^1\times\mathbb{Q}^2,~
\theta_t:=\theta_t^1\times\theta_t^2,
\de
and then $(\Omega, \sF, \mathbb{Q}, (\theta_t)_{t\in\mR})$ is a metric dynamical system, which will be used in the following.

\subsection{A random dynamical system}

Consider the system (\ref{slfasy}), i.e.
\ce\left\{\begin{array}{l}
\dot{u}^\e=\frac{1}{\e}A{u}^\e+\frac{1}{\e}U(u^\e,v^\e)+\frac{\sigma_1}{\e^{1/\alpha}}\dot{L}^{\alpha\pm}, \quad {u}_0^\e=u_0\in\mR^n,\\
\dot{v}^\e=B{v}^\e+V(u^\e, v^\e)+\sigma_2\dot{L}^\pm, \quad\qquad {v}_0^\e=v_0\in\mR^m.
\end{array}
\right.
\de

We make the following hypotheses:

\begin{enumerate}[($\bf{H_1}$)] 
\item There exists a $\gamma_1>0$ such that for any $x\in\mR^n$,
$$
\<Ax,x\>\leq -\gamma_1|x|^2.
$$
\end{enumerate}

\begin{enumerate}[($\bf{H_2}$)] 
\item There exist $\gamma_2, \gamma_3\geq0$ such that
\be
&&\|e^{Bt}\|\leq e^{-\gamma_2t}, \quad t\leq0,
\label{opecon1}\\
&&\|e^{Bt}\|\leq e^{-\gamma_3 t},  \quad t>0.
\label{opecon2}
\ee
\end{enumerate}
\begin{enumerate}[($\bf{H_3}$)] 
\item There exists a positive constant $L$ such that for all $(x_1, y_1), (x_2, y_2)\in \mR^n\times\mR^m$
\ce
|U(x_1, y_1)-U(x_2, y_2)|\leq L(|x_1-x_2|+|y_1-y_2|),
\de
and
\ce
|V(x_1, y_1)-V(x_2, y_2)|\leq L(|x_1-x_2|+|y_1-y_2|).
\de
\end{enumerate}
\begin{enumerate}[($\bf{H_4}$)] 
\item
$$
\gamma_1>L.
$$
\end{enumerate}
\begin{enumerate}[($\bf{H_5}$)] 
\item
\ce
&&\sup\limits_{(x, y)\in \mR^n\times\mR^m}|U(x,y)|=M_U,\\
&&\sup\limits_{(x, y)\in \mR^n\times\mR^m}|V(x,y)|=M_V.
\de
\end{enumerate}

Under the  assumptions $(\bf{H_1})$-$(\bf{H_3})$, the system (\ref{slfasy}) has a global unique
solution denoted by $(u^\e(t), v^\e(t))$(\cite{iw}). Define the solution operator
$\Psi^\e_t(u_0, v_0):=(u^\e(t), v^\e(t))$ for $t\geq0$, and then we know that $\Psi^\e$ is a random dynamical system.

\subsection{Random invariant manifolds}\label{rimep}

Introduce two auxiliary systems
\ce
&&\dif\zeta_t^\e=\frac{1}{\e}A\zeta_t^\e\dif t+\frac{\sigma_1}{\e^{1/\a}}\dif L^{\alpha\pm}_t, \quad \zeta^\e_0=u\in\mR^n,\\
&&\dif\varsigma_t=B\varsigma_t\dif t+\sigma_2\dif L^\pm_t, \quad\qquad \varsigma_0=v\in\mR^m.
\de
So, by \cite[Lemma 1]{zqd} and \cite[Example 3.7]{q1}, there exist two random vectors $\zeta^\e, \varsigma$ such that they are stationary solutions of two above equations. Set
\ce
&&\bar{u}^\e:=u^\e-\zeta^\e(\theta^1_{\cdot}\omega_1),\\
&&\bar{v}^\e:=v^\e-\varsigma(\theta^2_{\cdot}\omega_2),
\de
and then $(\bar{u}^\e, \bar{v}^\e)$ satisfy the following system
\be\left\{\begin{array}{l}
\dot{\bar{u}}^\e=\frac{1}{\e}A\bar{u}^\e+\frac{1}{\e}U(\bar{u}^\e+\zeta^\e(\theta^1_{\cdot}\omega_1),\bar{v}^\e+\varsigma(\theta^2_{\cdot}\omega_2)), \quad \bar{u}_0^\e=\bar{u}_0,
\\
\dot{\bar{v}}^\e=B{\bar{v}}^\e+V(\bar{u}^\e+\zeta^\e(\theta^1_{\cdot}\omega_1),\bar{v}^\e+\varsigma(\theta^2_{\cdot}\omega_2)), \quad\quad \bar{v}_0^\e=\bar{v}_0.
\end{array}
\right.
\label{slfasy2}
\ee
Moreover, $(\bar{u}_t^\e, \bar{v}_t^\e)$ generates a random dynamical system denoted by
$\bar{\Psi}_t^\e$ for $t\geq0$. 

Next, we construct a random invariant manifold with respect to $\bar{\Psi}^{\varepsilon}$. The method comes from \cite{q3}. We start with a key lemma.

\bl\label{equi}
Suppose that $(\bf{H_1})$-$(\bf{H_5})$ are satisfied. Then for $(\bar{u}_0,\bar{v}_0)\in\mR^n\times\mR^m$, there exists a $\e_0>0$ such that for $0<\e\leq\e_0$, the following integral equation has a unique solution $(\hat{u}_t^\e, \hat{v}_t^\e)$
\be
\left(\begin{array}{cccc}
\hat{u}_t^\e\\
\hat{v}_t^\e
\end{array}\right)=\left(\begin{array}{cccc}\int_{-\infty}^te^{\frac{A}{\e}(t-r)}\frac{1}{\e}U(\hat{u}_r^\e+\zeta^\e(\theta^1_r\omega_1),\hat{v}_r^\e+\varsigma(\theta^2_r\omega_2))dr\\
e^{Bt}\bar{v}_0-\int_t^0 e^{B(t-r)}V(\hat{u}_r^\e+\zeta^\e(\theta^1_r\omega_1),\hat{v}_r^\e+\varsigma(\theta^2_r\omega_2))dr\end{array}\right), t\leq 0,
\label{ineq}
\ee
\ce
\hat{v}_0^\e=\bar{v}_0.
\de
\el
\begin{proof}
Firstly, we introduce two following spaces
\ce
&&\cC^{-}_{\frac{\mu}{\e}}(\mR^n):=\left\{\phi\in\cC((-\infty,0], \mR^n): \sup\limits_{t\leq 0}e^{\frac{\mu}{\e} t}|\phi(t)|<\infty\right\},\\
&&\cC^{-}_{\frac{\mu}{\e}}(\mR^m):=\left\{\phi\in\cC((-\infty,0], \mR^m): \sup\limits_{t\leq 0}e^{\frac{\mu}{\e} t}|\phi(t)|<\infty\right\},
\de
where $\mu>0$ is a constant and $\gamma_1-\mu>L$. Let $\cC^-_{\frac{\mu}{\e}}(\mR^n\times\mR^m):=\cC^{-}_{\frac{\mu}{\e}}(\mR^n)\times\cC^{-}_{\frac{\mu}{\e}}(\mR^m)$ with the norm $\|z\|_{\cC^-_{\frac{\mu}{\e}}(\mR^n\times\mR^m)}=\sup\limits_{t\leq 0}e^{\frac{\mu}{\e} t}|z(t)|=\sup\limits_{t\leq 0}e^{\frac{\mu}{\e} t}(|u(t)|+|v(t)|)$ for $z=(u,v)\in\cC^-_{\frac{\mu}{\e}}(\mR^n\times\mR^m)$.

Set for $\hat{z}^\e=(\hat{u}^\e,\hat{v}^\e)\in\cC^-_{\frac{\mu}{\e}}(\mR^n\times\mR^m)$
\be
\cI(\hat{z}^\e)(t):=\left(\begin{array}{cccc}
\cI_1(\hat{z}^\e)(t)\\
\cI_2(\hat{z}^\e)(t)
\end{array}\right):=\left(\begin{array}{cccc}\int_{-\infty}^te^{\frac{A}{\e}(t-r)}\frac{1}{\e}U(\hat{u}_r^\e+\zeta^\e(\theta^1_r\omega_1),\hat{v}_r^\e+\varsigma(\theta^2_r\omega_2))dr\\
e^{Bt}\bar{v}_0-\int_t^0 e^{B(t-r)}V(\hat{u}_r^\e+\zeta^\e(\theta^1_r\omega_1),\hat{v}_r^\e+\varsigma(\theta^2_r\omega_2))dr\end{array}\right), 
\label{kmapp}
\ee
and then $\cI$ is well defined on $\cC^-_{\frac{\mu}{\e}}(\mR^n\times\mR^m)$. In fact, by $(\bf{H_1})$ $(\bf{H_2})$ $(\bf{H_5})$ it holds that for $\hat{z}^\e=(\hat{u}^\e,\hat{v}^\e)\in\cC^-_{\frac{\mu}{\e}}(\mR^n\times\mR^m)$,
\ce
\sup\limits_{t\leq 0}e^{\frac{\mu}{\e}t}|\cI_1(\hat{z}^\e)(t)|&\leq&\sup\limits_{t\leq 0}e^{\frac{\mu}{\e}t}\left|\int_{-\infty}^t e^{\frac{A}{\e}(t-r)}\frac{1}{\e}U(\hat{u}_r^\e+\zeta^\e(\theta^1_r\omega_1),\hat{v}_r^\e+\varsigma(\theta^2_r\omega_2))dr\right|\no\\
&\leq&\frac{M_U}{\e}\sup\limits_{t\leq 0}e^{\frac{\mu}{\e}t}\int_{-\infty}^te^{-\frac{\gamma_1}{\e}(t-r)}dr\leq\frac{M_U}{\gamma_1},
\de
and
\ce
\sup\limits_{t\leq 0}e^{\frac{\mu}{\e}t}|\cI_2(\hat{z}^\e)(t)|&\leq&\sup\limits_{t\leq 0}e^{\frac{\mu}{\e}t}|e^{Bt}\bar{v}_0|+\sup\limits_{t\leq 0}e^{\frac{\mu}{\e}t}\left|\int_t^0 e^{B(t-r)}V(\hat{u}_r^\e+\zeta^\e(\theta^1_r\omega_1),\hat{v}_r^\e+\varsigma(\theta^2_r\omega_2))dr\right|\\
&\leq&\sup\limits_{t\leq 0}e^{\frac{\mu}{\e}t}e^{-\gamma_2 t}|\bar{v}_0|+M_V\sup\limits_{t\leq 0}e^{\frac{\mu}{\e}t}\int_t^0 e^{-\gamma_2 (t-r)}dr\\
&\leq&|\bar{v}_0|+\frac{M_V}{\gamma_2}.
\de

In the following, we show that $\cI$ is a contractive mapping. For $\hat{z}^{\e,1}, \hat{z}^{\e,2}\in\cC^-_{\frac{\mu}{\e}}(\mR^n\times\mR^m)$, by $(\bf{H_3})$, one can obtain that 
\ce
\sup\limits_{t\leq 0}e^{\frac{\mu}{\e}t}\left|\cI_1(\hat{z}^{\e,1})(t)-\cI_1(\hat{z}^{\e,2})(t)\right|&\leq&\sup\limits_{t\leq 0}e^{\frac{\mu}{\e}t}\int_{-\infty}^te^{-\frac{\gamma_1}{\e}(t-r)}\frac{L}{\e}(|\hat{u}_r^{\e,1}-\hat{u}_r^{\e,2}|+|\hat{v}_r^{\e,1}-\hat{v}_r^{\e,2}|)\dif r\\
&\leq&\frac{L}{\gamma_1-\mu}\left(\sup\limits_{t\leq 0}e^{\frac{\mu}{\e}t}|\hat{z}_t^{\e,1}-\hat{z}_t^{\e,2}|\right),
\de
\ce
\sup\limits_{t\leq 0}e^{\frac{\mu}{\e}t}\left|\cI_2(\hat{z}^{\e,1})(t)-\cI_2(\hat{z}^{\e,2})(t)\right|&\leq&
\sup\limits_{t\leq 0}e^{\frac{\mu}{\e}t}\int_t^0 e^{-\gamma_2 (t-r)}L(|\hat{u}_r^{\e,1}-\hat{u}_r^{\e,2}|+|\hat{v}_r^{\e,1}-\hat{v}_r^{\e,2}|)\dif r\\
&\leq&\frac{\e L}{\mu-\e\gamma_2}\left(\sup\limits_{t\leq 0}e^{\frac{\mu}{\e}t}|\hat{z}_t^{\e,1}-\hat{z}_t^{\e,2}|\right).
\de
So, it holds that
\ce
\sup\limits_{t\leq 0}e^{\frac{\mu}{\e}t}\left|\cI(\hat{z}^{\e,1})(t)-\cI(\hat{z}^{\e,2})(t)\right|&\leq&\left(\frac{L}{\gamma_1-\mu}+\frac{\e L}{\mu-\e\gamma_2}\right)\left(\sup\limits_{t\leq 0}e^{\frac{\mu}{\e}t}|\hat{z}_t^{\e,1}-\hat{z}_t^{\e,2}|\right).
\de
Note that $\gamma_1-\mu> L$ and then 
$$
\frac{L}{\gamma_1-\mu}<1.
$$
Thus, there exists a $\e_0>0$ such that for any $0<\e\leq\e_0$,
$$
\frac{L}{\gamma_1-\mu}+\frac{\e L}{\mu-\e\gamma_2}<1,
$$
that is, $\cI$ is contractive. So, Eq.(\ref{ineq}) has a unique solution denoted as $(\hat{u}^\e,\hat{v}^\e)$. The proof is complete.
\end{proof}

Next, for $l\in(-\infty, 0]$, we rewrite Eq.(\ref{ineq}) as 
\ce
\left(\begin{array}{cccc}
\hat{u}_t^\e\\
\hat{v}_t^\e
\end{array}\right)=\left(\begin{array}{cccc}e^{\frac{A}{\e}(t-l)}\hat{u}_l^\e+\int_l^te^{\frac{A}{\e}(t-r)}\frac{1}{\e}U(\hat{u}_r^\e+\zeta^\e(\theta^1_r\omega_1),\hat{v}_r^\e+\varsigma(\theta^2_r\omega_2))dr\\
e^{B(t-l)}\hat{v}_l^\e+\int_l^te^{B(t-r)}V(\hat{u}_r^\e+\zeta^\e(\theta^1_r\omega_1),\hat{v}_r^\e+\varsigma(\theta^2_r\omega_2))dr
\end{array}\right), \quad l\leq t\leq 0.
\de
That is, the dynamic of Eq.(\ref{ineq}) is the same as that of the system (\ref{slfasy2}). 
Thus, it holds that for $t\leq 0$
\be
\hat{u}_t^\e(\t_s\omega)=\hat{u}_{t+s}^\e(\omega), \quad \hat{v}_t^\e(\t_s\omega)=\hat{v}_{t+s}^\e(\omega), \quad t+s\leq 0. 
\label{ineq2}
\ee

In the following, set
\ce
F^{\varepsilon}(\omega,\bar{v}_0)=\frac{1}{\e}\int_{-\infty}^0e^{-\frac{A}{\e}s}U\big(\hat{u}_s^\e+\zeta^\e(\theta^1_s\omega_1),\hat{v}_s^\e+\varsigma(\theta^2_s\omega_2)\big)ds,
\de
and then we study the properties of $F^{\varepsilon}(\omega,\bar{v}_0)$. First, by $(\bf{H_1})$ $(\bf{H_5})$ it holds that for $\bar{v}_0\in\mR^m$
\be
|F^{\varepsilon}(\omega,\bar{v}_0)|\leq\frac{M_U}{\e}\int_{-\infty}^0e^{\frac{\gamma_1}{\e}s}ds=\frac{M_U}{\gamma_1}.
\label{bouf}
\ee
Second, by the proof of Lemma \ref{equi} one can obtain that for $0<\varepsilon\leq\e_0$ and $\bar{v}^1_0, \bar{v}^2_0\in\mR^m$
$$
|F^{\varepsilon}(\omega,\bar{v}^1_0)-F^{\varepsilon}(\omega,\bar{v}^2_0)|\leq\frac{L}{(\gamma_1-\mu)\left[1-(\frac{L}{\gamma_1-\mu}+\frac{\e L}{\mu-\e\gamma_2})\right]}|\bar{v}^1_0-\bar{v}^2_0|.
$$
Third, it follows from (\ref{ineq2}) that for $t\leq0$,
\ce
F^\e(\t_t\omega,\bar{v}_0)=\int_{-\infty}^te^{\frac{A}{\e}(t-r)}\frac{1}{\e}U\big(\hat{u}_r^\e+\zeta^\e(\theta^1_r\omega_1),\hat{v}_r^\e+\varsigma(\theta^2_r\omega_2)\big)dr=\hat{u}_t^\e.
\de
Thus, define
$$
\bar{\mathcal{M}}^{\varepsilon}(\omega):=\{(F^{\varepsilon}(\omega, y),y), y\in\mR^m\},
$$
and then based on the property of $F^{\varepsilon}$, one can justify that $\bar{\mathcal{M}}^{\varepsilon}$ is a Lipschitz random invariant manifold with respect to $\bar{\Psi}^{\varepsilon}$. Moreover, we have the following result.

\bt\label{reesor}
Under the hypotheses $(\mathbf{H_1})$-$(\mathbf{H_5})$, it holds that for sufficiently small $\varepsilon$ and $(\bar{u}_0,\bar{v}_0)\in\mR^n\times\mR^m$, there exists $(\bar{\bar{u}}_0,\bar{\bar{v}}_0)\in\bar{\mathcal{M}}^{\varepsilon}$ such that 
\be
|\bar{\Psi}^{\varepsilon}_t (\bar{u}_0,\bar{v}_0)-\bar{\Psi}^{\varepsilon}_t(\bar{\bar{u}}_0,\bar{\bar{v}}_0)|\leq \frac{e^{-\frac{\mu}{\e} t}}{1-\(\frac{L}{\gamma_1-\mu}+\frac{\e L}{\mu-\e\gamma_2}\)}\(2(|\bar{u}_0|+|\bar{v}_0|)+2\frac{M_U}{\gamma_1}+\frac{M_V}{\gamma_2}\), \quad t\geq 0.
\label{esti}
\ee
\et
\begin{proof}
To prove the theorem, we need the following spaces:
\ce
&&\cC_{\frac{\mu}{\e}}(\mR^n):=\left\{\phi\in\cC(\mR, \mR^n): \sup\limits_{t\in\mR}e^{\frac{\mu}{\e} t}|\phi(t)|<\infty\right\},\\
&&\cC_{\frac{\mu}{\e}}(\mR^m):=\left\{\phi\in\cC(\mR, \mR^m): \sup\limits_{t\in\mR}e^{\frac{\mu}{\e} t}|\phi(t)|<\infty\right\},
\de
and $\cC_{\frac{\mu}{\e}}(\mR^n\times\mR^m):=\cC_{\frac{\mu}{\e}}(\mR^n)\times\cC_{\frac{\mu}{\e}}(\mR^m)$ with the norm $\|z\|_{\cC_{\frac{\mu}{\e}}(\mR^n\times\mR^m)}=\sup\limits_{t\in\mR}e^{\frac{\mu}{\e} t}|z(t)|=\sup\limits_{t\in\mR}e^{\frac{\mu}{\e} t}(|u(t)|+|v(t)|)$ for $z=(u,v)\in\cC_{\frac{\mu}{\e}}(\mR^n\times\mR^m)$.

Next, note that by the assumption $(\bf{H_1})$, the operator $I-|t|A$ is invertible. Thus, set 
\ce
\bar{z}_t^\e:=\left(\begin{array}{cccc}
\bar{u}_t^\e\\
\bar{v}_t^\e
\end{array}\right):=\left\{\begin{array}{cccc}
\((I-|t|A)^{-1}\bar{u}_0, \bar{v}_0\), \qquad t\leq0,\\
\bar{\Psi}^\e(t,\omega)(\bar{u}_0, \bar{v}_0), \qquad\quad\quad\quad t>0,
\end{array}\right.
\de
and
\ce
\bar{Z}_0(t):=\left\{\begin{array}{cccc}
-\bar{z}_t^\e+\cI(\bar{z}^\e)(t), \qquad\qquad\qquad\qquad\qquad t\leq0,\\
-(e^{\frac{A}{\e}t}\bar{u}_0,e^{B t}\bar{v}_0)+(e^{\frac{A}{\e}t}\cI_1(\bar{z}^\e)(0),e^{B t}\bar{v}_0), \quad t>0,
\end{array}\right.
\de
where $\cI$ is defined in (\ref{kmapp}). And then we consider the following integral equation
\be
\left(\begin{array}{cccc}
X_t^\e\\
Y_t^\e
\end{array}\right)
=\bar{Z}_0(t)+\left(\begin{array}{cccc}\int_{-\infty}^t e^{\frac{A}{\e}(t-r)}\frac{1}{\e}\[U(\bar{u}_r^\e+X_r^\e, \bar{v}_r^\e+Y_r^\e)-U(\bar{u}_r^\e, \bar{v}_r^\e)\]dr\\
-\int_t^\infty e^{B(t-r)}\[V(\bar{u}_r^\e+X_r^\e, \bar{v}_r^\e+Y_r^\e)-V(\bar{u}_r^\e, \bar{v}_r^\e)\]dr\end{array}\right), t\in\mR.
\label{ineq3}
\ee
For $Z^\e:=(X^\e, Y^\e)\in\cC_{\frac{\mu}{\e}}(\mR^n\times\mR^m)$, set
\ce
\cJ(Z^\e)(t)&&:=
\left(\begin{array}{cccc}
\cJ_1(Z^\e)(t)\\
\cJ_2(Z^\e)(t)
\end{array}\right)\\
&&:=\bar{Z}_0(t)+\left(\begin{array}{cccc}
\int_{-\infty}^t e^{\frac{A}{\e}(t-r)}\frac{1}{\e}\[U(\bar{u}_r^\e+X_r^\e, \bar{v}_r^\e+Y_r^\e)-U(\bar{u}_r^\e, \bar{v}_r^\e)\]dr\\
-\int_t^\infty e^{B(t-r)}\[V(\bar{u}_r^\e+X_r^\e, \bar{v}_r^\e+Y_r^\e)-V(\bar{u}_r^\e, \bar{v}_r^\e)\]dr\end{array}\right),
\de
and then $\cJ: \cC_{\frac{\mu}{\e}}(\mR^n\times\mR^m)\mapsto \cC_{\frac{\mu}{\e}}(\mR^n\times\mR^m)$ is well defined. In fact, for $Z^\e=(X^\e, Y^\e)\in\cC_{\frac{\mu}{\e}}(\mR^n\times\mR^m)$, by $(\bf{H_1})$-$(\bf{H_5})$ and (\ref{bouf}) we compute
\be
\sup\limits_{t\in\mR}e^{\frac{\mu}{\e} t}|\bar{Z}_0(t)|&\leq&\sup\limits_{t\leq0}e^{\frac{\mu}{\e} t}|\bar{Z}_0(t)|+\sup\limits_{t>0}e^{\frac{\mu}{\e} t}|\bar{Z}_0(t)|\no\\
&\leq&\sup\limits_{t\leq0}e^{\frac{\mu}{\e} t}(|\bar{z}_t^\e|+|\cI(\bar{z}^\e)(t)|)+\sup\limits_{t>0}e^{\frac{\mu}{\e} t}e^{-\frac{\gamma_1}{\e}t}\(|\bar{u}_0|+|\cI_1(\bar{z}^\e)(0)|\)\no\\
&\leq&\sup\limits_{t\leq0}e^{\frac{\mu}{\e} t}|\bar{z}_t^\e|+\sup\limits_{t\leq0}e^{\frac{\mu}{\e} t}|\cI(\bar{z}^\e)(t)|+|\bar{u}_0|+|\cI_1(\bar{z}^\e)(0)|\no\\
&\leq&|\bar{u}_0|+|\bar{v}_0|+|\bar{v}_0|+\frac{M_U}{\gamma_1}+\frac{M_V}{\gamma_2}+|\bar{u}_0|+\frac{M_U}{\gamma_1}\no\\
&=&2(|\bar{u}_0|+|\bar{v}_0|)+2\frac{M_U}{\gamma_1}+\frac{M_V}{\gamma_2}.
\label{esti0}
\ee
Besides, it follows from $(\bf{H_1})$-$(\bf{H_4})$ that
\be
&&\frac{1}{\e}\sup\limits_{t\in\mR}e^{\frac{\mu}{\e} t}\int_{-\infty}^t e^{-\frac{\gamma_1}{\e}(t-r)}|U(\bar{u}_r^\e+X_r^\e, \bar{v}_r^\e+Y_r^\e)-U(\bar{u}_r^\e, \bar{v}_r^\e)|dr\no\\
&\leq&\frac{L}{\e}\(\sup\limits_{t\in\mR}e^{\frac{\mu}{\e} t}|Z_t^\e|\)\sup\limits_{t\in\mR}\int_{-\infty}^t e^{(\frac{\mu}{\e}-\frac{\gamma_1}{\e})(t-r)}dr\no\\
&\leq&\frac{L}{\gamma_1-\mu}\(\sup\limits_{t\in\mR}e^{\frac{\mu}{\e} t}|Z_t^\e|\),
\label{esti1}
\ee
and
\be
&&\sup\limits_{t\in\mR}e^{\frac{\mu}{\e} t}\int_t^\infty e^{-\gamma_2(t-r)}|V(\bar{u}_r^\e+X_r^\e, \bar{v}_r^\e+Y_r^\e)-V(\bar{u}_r^\e, \bar{v}_r^\e)|dr\no\\
&\leq&L\(\sup\limits_{t\in\mR}e^{\frac{\mu}{\e} t}|Z_t^\e|\)\sup\limits_{t\in\mR}\int_t^\infty e^{(\frac{\mu}{\e}-\gamma_2)(t-r)}dr\no\\
&\leq&\frac{\e L}{\mu-\e\gamma_2}\(\sup\limits_{t\in\mR}e^{\frac{\mu}{\e} t}|Z_t^\e|\).
\label{esti2}
\ee
Thus, by combining (\ref{esti0}) (\ref{esti1}) with (\ref{esti2}), one can obtain that
\ce
\sup\limits_{t\in\mR}e^{\frac{\mu}{\e} t}|\cJ(Z^\e)(t)|\leq\sup\limits_{t\in\mR}e^{\frac{\mu}{\e} t}|\cJ_1(Z^\e)(t)|+\sup\limits_{t\in\mR}e^{\frac{\mu}{\e} t}|\cJ_2(Z^\e)(t)|<\infty.
\de

Next, for $Z^{\e,1}, Z^{\e,2}\in\cC_{\frac{\mu}{\e}}(\mR^n\times\mR^m)$, by the similar deduction to (\ref{esti1}) (\ref{esti2}) it holds that
\ce
\sup\limits_{t\in\mR}e^{\frac{\mu}{\e} t}|\cJ_1(Z^{\e,1})(t)-\cJ_1(Z^{\e,2})(t)|&\leq&\frac{L}{\gamma_1-\mu}\(\sup\limits_{t\in\mR}e^{\frac{\mu}{\e} t}|Z_t^{\e,1}-Z_t^{\e,2}|\),\\
\sup\limits_{t\in\mR}e^{\frac{\mu}{\e} t}|\cJ_2(Z^{\e,1})(t)-\cJ_2(Z^{\e,2})(t)|&\leq&
\frac{\e L}{\mu-\e\gamma_2}\(\sup\limits_{t\in\mR}e^{\frac{\mu}{\e} t}|Z_t^{\e,1}-Z_t^{\e,2}|\).
\de
Thus, we have that
\ce
 \sup\limits_{t\in\mR}e^{\frac{\mu}{\e} t}|\cJ(Z^{\e,1})(t)-\cJ(Z^{\e,2})(t)|
 \leq\(\frac{L}{\gamma_1-\mu}+\frac{\e L}{\mu-\e\gamma_2}\)\(\sup\limits_{t\in\mR}e^{\frac{\mu}{\e} t}|Z_t^{\e,1}-Z_t^{\e,2}|\).
\de
So, for $0<\e\leq\e_0$, $\cJ: \cC_{\frac{\mu}{\e}}(\mR^n\times\mR^m)\mapsto \cC_{\frac{\mu}{\e}}(\mR^n\times\mR^m)$ is contractive. That is, Eq.(\ref{ineq3}) has a unique solution denoted as $Z^\e=(X^\e, Y^\e)$. Moreover, 
\ce
\sup\limits_{t\in\mR}e^{\frac{\mu}{\e} t}|Z_t^\e|\leq \frac{1}{1-\(\frac{L}{\gamma_1-\mu}+\frac{\e L}{\mu-\e\gamma_2}\)}\(2(|\bar{u}_0|+|\bar{v}_0|)+2\frac{M_U}{\gamma_1}+\frac{M_V}{\gamma_2}\)
\de
and then 
\be
|Z_t^\e|\leq \frac{e^{-\frac{\mu}{\e} t}}{1-\(\frac{L}{\gamma_1-\mu}+\frac{\e L}{\mu-\e\gamma_2}\)}\(2(|\bar{u}_0|+|\bar{v}_0|)+2\frac{M_U}{\gamma_1}+\frac{M_V}{\gamma_2}\)
, \quad t\geq 0.
\label{esti3}
\ee

Set
\ce
\bar{\bar{z}}_t^\e:=\left(\begin{array}{cccc}
\bar{\bar{u}}_t^\e\\
\bar{\bar{v}}_t^\e\end{array}\right):=\left(\begin{array}{cccc}
\bar{u}_t^\e\\
\bar{v}_t^\e\end{array}\right)+\left(\begin{array}{cccc}
X_t^\e\\
Y_t^\e\end{array}\right),
\de
and then by simple calculation, it holds that $\bar{\bar{z}}_t^\e=(\bar{\bar{u}}_t^\e,\bar{\bar{v}}_t^\e)$ solves uniquely the following equation
\ce
\bar{\bar{z}}_t^\e=\left\{\begin{array}{cccc}\cI(\bar{\bar{z}}^\e)(t),\qquad\qquad t\leq0,\\
\bar{\Psi}^\e(t,\omega)(\bar{\bar{u}}_0^\e, \bar{\bar{v}}_0^\e), \quad t>0.
\end{array}\right.
\de
Thus, by Lemma \ref{equi}, we know that $\bar{\bar{z}}_t^\e$ solves Eq.(\ref{ineq}) for $t\leq0$. In particular, $\bar{\bar{u}}_0^\e=F^{\e}(\omega, \bar{\bar{v}}_0^\e)$, which yields that $(\bar{\bar{u}}_0^\e, \bar{\bar{v}}_0^\e)\in\bar{\cM}^\e(\omega)$. So, one can take $\bar{\bar{z}}_0=(\bar{\bar{u}}_0^\e, \bar{\bar{v}}_0^\e)$. Since $\bar{\bar{z}}_t^\e=\bar{\Psi}^\e(t,\omega)(\bar{\bar{u}}_0^\e, \bar{\bar{v}}_0^\e)$ for $t>0$, $\bar{\bar{z}}_t^\e=\bar{\Psi}^\e(t,\omega)\bar{\bar{z}}_0$ for $t>0$.  Note that $\bar{\bar{u}}_t^\e-\bar{u}_t^\e=X_t^\e, \bar{\bar{v}}_t^\e-\bar{v}_t^\e=Y_t^\e$. Thus, by (\ref{esti3}), it holds that
\ce
|\bar{\Psi}^\e(t,\omega)\bar{z}_0-\bar{\Psi}^\e(t,\omega)\bar{\bar{z}}_0|&=&|(\bar{u}_t^\e, \bar{v}_t^\e)-(\bar{\bar{u}}_t^\e, \bar{\bar{v}}_t^\e)|=|Z_t^\e|\\
&\leq&\frac{e^{-\frac{\mu}{\e} t}}{1-\(\frac{L}{\gamma_1-\mu}+\frac{\e L}{\mu-\e\gamma_2}\)}\(2(|\bar{u}_0|+|\bar{v}_0|)+2\frac{M_U}{\gamma_1}+\frac{M_V}{\gamma_2}\), \quad t\geq 0.
\de
The proof is complete.
\end{proof}

Based on the relationship between \eqref{slfasy} and \eqref{slfasy2}, it holds that the system \eqref{slfasy} has a random invariant manifold
\ce
\mathcal{M}^{\varepsilon}(\omega)=\{(F^{\varepsilon}(\omega, y)+\zeta^{\varepsilon}(\omega_1),y+\varsigma(\omega_2)), y\in\mR^m\}.
\de

\subsection{A reduction system on the random invariant manifold $\mathcal{M}^{\varepsilon}$}

By Theorem \ref{reesor}, we can get the following reduced system approximating the system \eqref{slfasy}.

\bt\label{resyth}
Suppose that these assumptions $(\mathbf{H_1})$-$(\mathbf{H_5})$ hold. Then for any solution $z^{\varepsilon}_t= (u^{\varepsilon}_t, v^{\varepsilon}_t)$ to the system \eqref{slfasy} with the initial data $z^{\varepsilon}_0= (u_0, v_0)$, there exists the following reduced low dimensional system on the random invariant manifold $\mathcal{M}^{\varepsilon}$
\begin{equation}
\label{redda}
\left\{
\begin{aligned}
{\tilde{u}^\varepsilon _t} &=F^{\varepsilon}(\theta_t \omega, {\tilde{v}^\varepsilon _t}-\varsigma(\theta^{2}_t\omega_2))+\zeta^{\varepsilon}(\theta^1_t\omega_1),\\
d{\tilde{v}^\varepsilon _t} &= B\tilde{v}^\varepsilon_t dt+{V}({\tilde{u}^\varepsilon_t},{\tilde{v}^\varepsilon _t})dt +\sigma_2 dL_t^\pm,
\end{aligned}
\right.
\end{equation}
such that for $0<\e\leq \e_0$, we have
\be
|z^{\varepsilon}(t,\omega)-\tilde{z}^{\varepsilon}(t,\omega)| \leq \frac{e^{-\frac{\mu}{\e} t}}{1-\(\frac{L}{\gamma_1-\mu}+\frac{\e L}{\mu-\e\gamma_2}\)}\(2(|u_0-u|+|v_0-v|)+2\frac{M_U}{\gamma_1}+\frac{M_V}{\gamma_2}\), t\geq 0,
\label{orrees}
\ee
where $ \tilde{z}^{\varepsilon}_t= (\tilde{u}^{\varepsilon}_t, \tilde{v}^{\varepsilon}_t)$  is the solution of the low dimensional system \eqref{redda} with the initial value $\tilde{z}^{\varepsilon}_0$.
\et

Note that in the system (\ref{redda}), ${\tilde{u}^\varepsilon _t}$ can be represented by ${\tilde{v}^\varepsilon _t}$. Thus, the system (\ref{redda}) is essentially decided by ${\tilde{v}^\varepsilon _t}$. That is, the dimension of the system (\ref{redda}) is $m$. However, the dimension of the system \eqref{slfasy} is $m+n$. Therefore, the system (\ref{redda}) is a reduced low dimensional system.

\br
By the estimate (\ref{orrees}), we know that when $\e$ is enough small or $t$ is sufficiently large, the system (\ref{redda}) will approximate the system \eqref{slfasy}.
\er

\br
Note that the estimate (\ref{orrees}) is different from the estimate in \cite[Theorem 3.2]{q4} and (3.17) in \cite{zqd}. In fact, the latter estimates omit the property of the reduced system. Besides, in the system (\ref{slfasy}), $L^\pm$ can be replaced by a $m$-dimensional two-sided symmetric $\a$-stable process. Moreover, by the similar deduction to that in Theorem \ref{resyth}, we can also obtain the estimate (\ref{orrees}).
\er

\section{An approximate filtering on the invariant manifold}\label{filter}

In the section we introduce nonlinear filtering problems for the system (\ref{slfasy}) and the reduced
system (\ref{redda}) on the random  invariant manifold, and then study their relation.

\subsection{Nonlinear filtering problems}\label{nonfil}

In the subsection we introduce nonlinear filtering problems for the system (\ref{slfasy}) and the reduced system (\ref{redda}).

For $T>0$, we take an observation system as follows
\ce
w^{\e}_t= W_t+\int_0^t H(u^\e_s, v^\e_s) \dif s, \quad t\in[0,T],
\de
where $W$ is a $l$-dimensional standard Brownian motion. Here $W$ may be either independent of $L^{\a\pm}$ and $L^\pm$, or correlated with $L^{\a\pm}$ and $L^\pm$. Assume:
\begin{enumerate}[($\bf{H_6}$)] 
\item $H$ is bounded and Lipschitz continuous in $(u,v)$ with the Lipschitz constant $\|H\|_{Lip}$.
\end{enumerate}

Set
$$
(\chi_t^\e)^{-1}:=\exp\left\{-\int_0^t H(u_s^\e, v_s^\e)\dif W_s-\frac12\int_0^t|H(u_s^\e, v_s^\e)|^2\dif
s\right\},
$$
and then we know that $(\chi_t^\e)^{-1}$ is an exponential martingale under $\mQ$. And define
$$
\dif \mQ^\e:=(\chi^\e_T)^{-1}\dif \mQ,
$$
and it holds that $\mQ^\e$ is a probability measure and $w^{\e}$ is a standard Brownian motion under $\mQ^\e$. Besides, set for $\Phi\in \cB(\mR^n\times\mR^m)$ (the set of
all real-valued uniformly bounded Borel-measurable functions on $\mR^n\times\mR^m$)
\ce
 \mQ_t^\e(\Phi) &:=&\mE^\e[\Phi(u_t^\e, v_t^\e)\chi^\e_t|\mathcal{W}_t^\e],\\
 \Pi_t^\e(\Phi) &:=& \mE[\Phi(u_t^\e, v_t^\e)|\mathcal{W}_t^\e],
 \de
 where  $\mE^\e, \mE$ stand for the expectation under $\mQ^\e$ and $\mQ$, respectively, $\mathcal{W}_t^\e \triangleq\sigma(w_s^\e:
 0\leq s \leq t) \vee \cN$ and $\cN$ is the collection of all $\mQ$-measure zero sets. Moreover, by the Kallianpur-Striebel formula it holds that
\ce
\Pi^{\e}_t(\Phi)=\frac{\mQ^{\e}_t(\Phi)}{\mQ^{\e}_t(1)}.
\de

Next, we study the nonlinear filtering problem for $(\tilde{u}^\e, \tilde{v}^\e)$. Set
\ce
\tilde{\chi}^\e_t:=\exp\left\{\int_0^t H(\tilde{u}^\e_s, \tilde{v}_s^\e)\dif
w^{\e}_s-\frac12\int_0^t|H(\tilde{u}^\e_s, \tilde{v}_s^\e)|^2\dif s\right\},
\de
and then $\tilde{\chi}^\e_t$ is an exponential martingale under $\mQ^\e$. Thus, we define the
nonnormalized filtering for $(\tilde{u}^\e_t, \widetilde{v}_t^\e)$ by
$$
\tilde{\mQ}_t^\e(\Phi) :=\mE^\e[\Phi(\tilde{u}^\e_t, \widetilde{v}_t^\e)\tilde{\chi}^\e_t|\mathcal{W}_t^\e].
$$
And set
$$
\tilde{\Pi}_t^\e(\Phi):=\frac{\tilde{\mQ}_t^\e(\Phi)}{\tilde{\mQ}_t^\e(1)},
$$
and then we will prove that $\tilde{\Pi}^\e_t$ could be understood as the nonlinear filtering problem for
$(\tilde{u}^\e_t, \tilde{v}_t^\e)$ with respect to $\mathcal{W}_t^\e$.

\subsection{The relation between $\Pi^{\e}_t$ and $\tilde{\Pi}^\e_t$}

In the subsection we will prove that $\tilde{\Pi}^\e_t$ converges to $\Pi^{\e}_t$ as $\e\rightarrow0$. Here let $\cC_b^1(\mR^n\times\mR^m)$ denote the collection of all functions which themselves and their first-order derivatives are uniformly bounded. We introduce the following norm for $\Phi\in\cC_b^1(\mR^n\times\mR^m)$:
\ce
\|\Phi\|_{\cC_b^1(\mR^n\times\mR^m)}=\max\limits_{(x,y)\in\mR^n\times\mR^m}|\Phi(x,y)|+\max\limits_{(x,y)\in\mR^n\times\mR^m}|\triangledown\Phi(x,y)|,
\de
where $\triangledown$ stands for the gradient operator. 

\bt\label{filcon} 
Under  $(\bf{H_1})$-$(\bf{H_6})$,  there exists a positive constant $C$ independent of $\e$ such that for $0<\e\leq\e_0$ and $\Phi\in \cC_b^1(\mR^n\times\mR^m)$,
\be
\mE|\Pi^{\e}_t(\Phi)- \tilde{\Pi}_t^\e(\Phi)|^p&\leq&\|\Phi\|_{\cC_b^1(\mR^n\times\mR^m)}^{p}\frac{C}{\[1-\left(\frac{L}{\gamma_1-\mu}+\frac{\e L}{\mu-\e\gamma_2}\right)\]^{p}}\(e^{-\frac{4p\mu}{\e} t}+\frac{\e}{4\mu p}\)^{1/4}, t\in[0,T], p>1.\no\\
\label{thes}
\ee
\et
\begin{proof}

For $\Phi\in \cC^1_b(\mR^n\times\mR^m)$, we compute that
\be
\mE|\Pi^{\e}_t(\Phi)-
\tilde{\Pi}_t^\e(\Phi)|^{p}&=&\mE\left|\frac{\mQ^{\e}_t(\Phi)-\tilde{\mQ}^{\e}_t(\Phi)}{\tilde{\mQ}^{\e}_t(1)}-\Pi^{\e}_t(\Phi)\frac{\mQ^{\e}_t(1)-\tilde{\mQ}^{\e}_t(1)}{\tilde{\mQ}^{\e}_t(1)}\right|^{p}\no\\
&\leq&2^{p-1}\mE\left|\frac{\mQ^{\e}_t(\Phi)-\tilde{\mQ}^{\e}_t(\Phi)}{\tilde{\mQ}^{\e}_t(1)}\right|^{p}+2^{p-1}\mE\left|\Pi^{\e}_t(\Phi)\frac{\mQ^{\e}_t(1)-\tilde{\mQ}^{\e}_t(1)}{\tilde{\mQ}^{\e}_t(1)}\right|^{p}\no\\
&\leq&2^{p-1}\left(\mE\left|\mQ^{\e}_t(\Phi)-\tilde{\mQ}^{\e}_t(\Phi)\right|^{2p}\right)^{1/2}\left(\mE\left|\tilde{\mQ}^{\e}_t(1)\right|^{-2p}\right)^{1/2}\no\\
&&+2^{p-1}\|\Phi\|_{\cC_b^1(\mR^n\times\mR^m)}^{p}\left(\mE\left|\mQ^{\e}_t(1)-\tilde{\mQ}^{\e}_t(1)\right|^{2p}\right)^{1/2}\left(\mE\left|\tilde{\mQ}^{\e}_t(1)\right|^{-2p}\right)^{1/2}.\no\\
\label{diffes}
\ee
In the following, it is the main task to estimate $\mE\left|\mQ^{\e}_t(\Phi)-\tilde{\mQ}^{\e}_t(\Phi)\right|^{2p}$ and $\mE\left|\tilde{\mQ}^{\e}_t(1)\right|^{-2p}$. 

By the similar deduction to that of \cite[Lemma 5.1]{q3}, it holds that 
\be
\mE\left|\tilde{\mQ}^{\e}_t(1)\right|^{-2p}\leq\exp\left\{(8p^2+2p+1)CT/2\right\}.\label{testi1}
\ee
Next, we are devoted to dealing with $\mE\left|\mQ^{\e}_t(\Phi)-\tilde{\mQ}^{\e}_t(\Phi)\right|^{2p}$. Following up the line in \cite[Lemma 5.2]{q3}, one can obtain that 
\be
\mE\left|\mQ^{\e}_t(\Phi)-\tilde{\mQ}^{\e}_t(\Phi)\right|^{2p}&\leq&\|\Phi\|_{\cC_b^1(\mR^n\times\mR^m)}^{2p}\frac{C}{\[1-\left(\frac{L}{\gamma_1-\mu}+\frac{\e L}{\mu-\e\gamma_2}\right)\]^{2p}}\no\\
&&\cdot\(\mE^\e\(2(|u_0-u|+|v_0-v|)+2\frac{M_U}{\gamma_1}+\frac{M_V}{\gamma_2}\)^{8p}\)^{1/4}\(e^{-\frac{4p\mu}{\e} t}+\frac{\e}{4\mu p}\)^{1/2}\no\\
&\leq&\|\Phi\|_{\cC_b^1(\mR^n\times\mR^m)}^{2p}\frac{C}{\[1-\left(\frac{L}{\gamma_1-\mu}+\frac{\e L}{\mu-\e\gamma_2}\right)\]^{2p}}\(e^{-\frac{4p\mu}{\e} t}+\frac{\e}{4\mu p}\)^{1/2}.
\label{testi2}
\ee
Thus, combining (\ref{diffes}) with (\ref{testi1})-(\ref{testi2}), we obtain (\ref{thes}). The proof is complete.
\end{proof}

\br
(\ref{thes}) indicates that when $\e$ goes to zero, $\tilde{\Pi}_t^\e$ approximates $\Pi^{\e}_t$. Therefore, $\tilde{\Pi}^\e$ could be understood as the nonlinear filtering problem for $(\tilde{u}^\e_t, \tilde{v}_t^\e)$ with respect to $\mathcal{W}_t^\e$.
\er

\br
If we take $\Phi(x,y)=\psi(y)$, Theorem \ref{filcon} is Theorem 3 in \cite{zqd}. Therefore, our result is more general.
\er

\section{The reduction for $\e=0$}\label{epze}

In the section, we observe the system (\ref{slfasy}) with $\sigma_2=0$, i.e.
\be\left\{\begin{array}{l}
\dot{u}^\e=\frac{1}{\e}A{u}^\e+\frac{1}{\e}U(u^\e,v^\e)+\frac{\sigma_1}{\e^{1/\alpha}}\dot{L}^{\alpha\pm}, \\
\dot{v}^\e=B{v}^\e+V(u^\e, v^\e).
\end{array}
\right.
\label{sig20}
\ee

\subsection{A reduced system for $\e=0$} 

In the subsection, we investigate the system (\ref{sig20}) for $\e=0$. To do this, we scale the time $t\rightarrow \e t$ and rewrite the system (\ref{sig20}) as
\be\left\{\begin{array}{l}
\dot{\check{u}}^\e=A{\check{u}}^\e+U(\check{u}^\e,\check{v}^\e)+\sigma_1\dot{\check{L}}^{\alpha\pm}, \quad {\check{u}}_0^\e=u_0\in\mR^n,\\
\dot{\check{v}}^\e=\e B{\check{v}}^\e+\e V(\check{u}^\e, \check{v}^\e), \quad\quad\qquad {\check{v}}_0^\e=v_0\in\mR^m,
\label{slfasysc}
\end{array}
\right.
\ee
where $\check{u}^\e_t:=u^\e_{\e t}, \check{v}^\e_t:=v^\e_{\e t}, \check{L}^{\alpha\pm}_t:=\frac{L^{\alpha\pm}_{\e t}}{\e^{1/\alpha}}$. By the scaling property of $\a$-stable processes, one can obtain that $\check{L}^{\alpha\pm}$ is still a two-sided $\a$-stable process. Introducing an auxiliary system
\ce
\dif\zeta_t=A\zeta_t\dif t+\sigma_1\dif \check{L}^{\alpha\pm}_t, \quad \zeta_0=u\in\mR^n,
\de
by \cite[Lemma 1]{zqd}, we know that there exists a random variable $\zeta$ such that it is a stationary solution of the above equation. Set
\ce
\bar{\check{u}}^\e:=\check{u}^\e-\zeta(\theta^1_{\cdot}\omega_1),
\de
and then $(\bar{\check{u}}^\e, \check{v}^\e)$ satisfy the following system
\be\left\{\begin{array}{l}
\dot{\bar{\check{u}}}^\e=A\bar{\check{u}}^\e+U(\bar{\check{u}}^\e+\zeta(\theta^1_{\cdot}\omega_1),\check{v}^\e), \quad \bar{\check{u}}_0^\e=\bar{u}_0\in\mR^n,
\\
\dot{\check{v}}^\e=\e B{\check{v}}^\e+\e V(\bar{\check{u}}^\e+\zeta(\theta^1_{\cdot}\omega_1),\check{v}^\e), \quad \check{v}_0^\e=v_0\in\mR^m.
\end{array}
\right.
\label{slfasy22}
\ee

Next, set
\ce
\check{F}^{\varepsilon}(\omega,v_0):=\int_{-\infty}^0e^{-As}U\big(\hat{\check{u}}_s^\e+\zeta(\theta^1_s\omega_1),\check{v}_s^\e\big)ds,
\de
where $\hat{\check{u}}_s^\e$ is the solution of an integral equation similar to Eq.(\ref{ineq}), and then by the similar deduction to that in Subsection \ref{rimep}, it holds that for $0<\varepsilon\leq\e_0$ and $v^1_0, v^2_0\in\mR^m$
$$
|\check{F}^{\varepsilon}(\omega,v^1_0)-\check{F}^{\varepsilon}(\omega,v^2_0)|\leq\frac{L}{(\gamma_1-\mu)\left[1-(\frac{L}{\gamma_1-\mu}+\frac{\e L}{\mu-\e\gamma_2})\right]}|v^1_0-v^2_0|.
$$
Again set 
$$
\check{\cM}^\e(\omega):=\left\{\big(\check{F}^{\varepsilon}(\omega,y)+\zeta(\omega_1), y\big), y\in\mR^m\right\},
$$
and then by the similar deduction to that in Theorem \ref{resyth}, we know that for any solution $\check{z}^{\varepsilon}_t= (\check{u}^{\varepsilon}_t, \check{v}^{\varepsilon}_t)$ to the system \eqref{slfasysc} with the initial value $\check{z}^{\varepsilon}_0= (u_0, v_0)$, there exists the following reduced low dimensional system on the random invariant manifold $\check{\cM}^\e$
\begin{equation}
\label{reddasc}
\left\{
\begin{aligned}
{\tilde{\check{u}}^\varepsilon _t} &=\check{F}^{\varepsilon}(\theta_t \omega, {\tilde{\check{v}}^\varepsilon _t})+\zeta(\theta^1_t\omega_1),\\
d{\tilde{\check{v}}^\varepsilon _t} &= \e B\tilde{\check{v}}^\varepsilon_t dt+\e{V}({\tilde{\check{u}}^\varepsilon_t},{\tilde{\check{v}}^\varepsilon _t})dt,
\end{aligned}
\right.
\end{equation}
such that for $0<\e\leq \e_0$, we have
\be
|\check{z}^{\varepsilon}(t,\omega)-\tilde{\check{z}}^{\varepsilon}(t,\omega)| \leq \frac{e^{-\mu t}}{1-\left(\frac{L}{\gamma_1-\mu}+\frac{\e L}{\mu-\e\gamma_2}\right)}\(2|u_0-u|+2\frac{M_U}{\gamma_1}+\frac{M_V}{\gamma_2}\), \quad t\geq 0,
\label{orreessc}
\ee
where $ \tilde{\check{z}}^{\varepsilon}_t= (\tilde{\check{u}}^{\varepsilon}_t, \tilde{\check{v}}^{\varepsilon}_t)$  is the solution of the low dimensional system \eqref{reddasc} with the initial value $\tilde{\check{z}}^{\varepsilon}_0$.

As $\e\rightarrow 0$, these systems (\ref{slfasysc}) and (\ref{slfasy22}) become
\be\left\{\begin{array}{l}
\dot{\check{u}}^0=A{\check{u}}^0+U(\check{u}^0,\check{v}^0)+\sigma_1\dot{\check{L}}^{\alpha\pm}, ~\quad\quad{\check{u}}_0^0=u_0,\\
\dot{\check{v}}^0=0,  \quad\qquad\qquad\qquad\qquad\qquad\qquad {\check{v}}_0^0=v_0,
\label{slfasysce0}
\end{array}
\right.
\ee
and
\be\left\{\begin{array}{l}
\dot{\bar{\check{u}}}^0=A\bar{\check{u}}^0+U(\bar{\check{u}}^0+\zeta(\theta^1_{\cdot}\omega_1),\check{v}^0), \qquad \bar{\check{u}}_0^0=\bar{u}_0\in\mR^n,
\\
\dot{\check{v}}^0=0, \quad\qquad\qquad\qquad\qquad\qquad\qquad \check{v}_0^0=v_0\in\mR^m,
\end{array}
\right.
\label{slfasy2e0}
\ee
respectively. Set
\ce
\check{F}^{0}(\omega,v_0):=\int_{-\infty}^0e^{-As}U\big(\hat{\check{u}}_s^0+\zeta(\theta^1_s\omega_1),v_0\big)ds,
\de
where $\hat{\check{u}}_s^0$ is the solution of an integral equation similar to Eq.(\ref{ineq}), and then by the similar deduction to that in Subsection \ref{rimep}, we have that for $v_0^1, 
v_0^2\in\mR^m$
\be
|\check{F}^0(\omega,v^1_0)-\check{F}^0(\omega,v^2_0)|\leq\frac{L}{\gamma_1-\mu-L}|v^1_0-v^2_0|.
\label{lipconsc0}
\ee
Put
$$
\check{\cM}^0(\omega):=\left\{\big(\check{F}^{0}(\omega,y)+\zeta(\omega_1), y\big), y\in\mR^m\right\},
$$
and by the same deduction to that in Theorem \ref{resyth}, one can obtain that there exists the following reduced low dimensional system on the random invariant manifold $\check{\cM}^0$
\begin{equation}
\label{reddasce0}
\left\{
\begin{aligned}
{\tilde{\check{u}}^0 _t} &=\check{F}^0(\theta_t \omega, {\tilde{\check{v}}^0 _t})+\zeta(\theta^1_t\omega_1),\\
d{\tilde{\check{v}}^0 _t} &= 0,
\end{aligned}
\right.
\end{equation}
such that 
\be
|\check{z}^0(t,\omega)-\tilde{\check{z}}^0(t,\omega)| \leq \frac{\gamma_1-\mu}{\gamma_1-\mu-L}e^{-\mu t}\(2|u_0-u|+2\frac{M_U}{\gamma_1}+\frac{M_V}{\gamma_2}\), \quad t\geq 0,
\label{orreessce0}
\ee
where $\check{z}^0_t= (\check{u}^0_t, \check{v}^0_t)$ and $\tilde{\check{z}}^0_t= (\tilde{\check{u}}^0_t, \tilde{\check{v}}^0_t)$  are the solutions of the systems \eqref{slfasysce0}  and \eqref{reddasce0} with the initial data $\check{z}^0_0=(u_0,v_0)$ and $\tilde{\check{z}}^0_0$, respectively.

\subsection{The relation between the system \eqref{slfasysc} and the reduced system \eqref{reddasce0}}

In the subsection we investigate the relation between the system \eqref{slfasysc} and the system \eqref{reddasce0}. Firstly, we need the following lemma.

\bl\label{madis}
Assume that $(\bf{H_1})$--$(\bf{H_5})$ hold. Then for $v_0\in\mR^m$
\ce
|\check{F}^{\varepsilon}(\omega, v_0)-\check{F}^0(\omega, v_0)|\leq C\beta(\e),
\de
where $C>0$ is a constant independent of $\e$ and 
\ce
\beta(\e):=e^{\mu t_0}\(\frac{1}{\gamma_1-\e \gamma_2}e^{-\e \gamma_2 t_0}-\frac{1}{\gamma_1}\)+\(\frac{1}{\gamma_1-\e \gamma_2}-\frac{1}{\gamma_1}\),\quad
t_0:=\frac{1}{\e \gamma_2}\log\frac{(\mu-\e \gamma_2)\gamma_1}{(\gamma_1-\e \gamma_2)\mu}.
\de
\el

Since the proof of the above lemma is similar to that of \cite[Theorem 5.1]{Fu}, we omit it.

\bt\label{mali}
Suppose that $(\bf{H_1})$-$(\bf{H_5})$ are satisfied. Then for $0<\e\leq\e_0$
\be
|\check{z}^{\varepsilon}_t-\tilde{\check{z}}^0_t|&\leq&\frac{e^{-\mu t}}{1-\left(\frac{L}{\gamma_1-\mu}+\frac{\e L}{\mu-\e\gamma_2}\right)}\(2|u_0-u|+2\frac{M_U}{\gamma_1}+\frac{M_V}{\gamma_2}\)\no\\
&&+\frac{1-\frac{\e L}{\mu-\e\gamma_2}}{1-\left(\frac{L}{\gamma_1-\mu}+\frac{\e L}{\mu-\e\gamma_2}\right)}\frac{|B v_0|+M_V}{\gamma_3}(1-e^{-\e\gamma_3 t})\no\\
&&+C\beta(\e), \qquad\qquad\qquad\qquad\qquad t\geq 0.
\label{malidi}
\ee
\et
\begin{proof}
By (\ref{orreessc}), it holds that
\be
|\check{z}^{\varepsilon}_t-\tilde{\check{z}}^0_t|&\leq&|\check{z}^{\varepsilon}_t-\tilde{\check{z}}^{\e}_t|+|\tilde{\check{z}}^{\varepsilon}_t-\tilde{\check{z}}^0_t|\no\\
&\leq&\frac{e^{-\mu t}}{1-\left(\frac{L}{\gamma_1-\mu}+\frac{\e L}{\mu-\e\gamma_2}\right)}\(2|u_0-u|+2\frac{M_U}{\gamma_1}+\frac{M_V}{\gamma_2}\)+|\tilde{\check{z}}^{\varepsilon}_t-\tilde{\check{z}}^0_t|.\no\\
\label{ches1}
\ee
Then we estimate $|\tilde{\check{z}}^{\varepsilon}_t-\tilde{\check{z}}^0_t|$. Note that $\tilde{\check{z}}^{\varepsilon}, \tilde{\check{z}}^0$ satisfy these systems (\ref{reddasc}) and (\ref{reddasce0}), respectively. 
Thus, it follows from Lemma \ref{madis} that
\be
|\tilde{\check{z}}^{\varepsilon}_t-\tilde{\check{z}}^0_t|&=&|\tilde{\check{u}}^{\varepsilon}_t-\tilde{\check{u}}^0_t|+|\tilde{\check{v}}^{\varepsilon}_t-\tilde{\check{v}}^0_t|=|\check{F}^{\varepsilon}(\theta_t \omega, {\tilde{\check{v}}^\varepsilon _t})-\check{F}^0(\theta_t \omega, {\tilde{\check{v}}^0 _t})|+|\tilde{\check{v}}^{\varepsilon}_t-\tilde{\check{v}}^0_t|\no\\
&\leq&|\check{F}^{\varepsilon}(\theta_t \omega, {\tilde{\check{v}}^\varepsilon _t})-\check{F}^{\e}(\theta_t \omega, v_0)|+|\check{F}^{\varepsilon}(\theta_t \omega, v_0)-\check{F}^0(\theta_t \omega, v_0)|+|\tilde{\check{v}}^{\varepsilon}_t-v_0|\no\\
&\leq&\frac{L}{(\gamma_1-\mu)\left[1-(\frac{L}{\gamma_1-\mu}+\frac{\e L}{\mu-\e\gamma_2})\right]}|\tilde{\check{v}}^{\varepsilon}_t-v_0|+C\beta(\e)+|\tilde{\check{v}}^{\varepsilon}_t-v_0|.
\label{ches2}
\ee
In the following, we are devoted to computing $|\tilde{\check{v}}^{\varepsilon}_t-v_0|$. Based on (\ref{reddasc}), it holds that 
\be
|\tilde{\check{v}}^{\varepsilon}_t-v_0|&\leq& |e^{\e B t}v_0-v_0|+\e\left|\int_0^te^{\e B(t-s)}{V}({\tilde{\check{u}}^\varepsilon_s},{\tilde{\check{v}}^\varepsilon _s})ds\right|\no\\
&\leq&\e|B v_0|\int_0^t\|e^{\e B s}\|ds+\e M_V\int_0^t\|e^{\e B(t-s)}\|ds\no\\
&\leq&\e|B v_0|\int_0^t e^{-\e\gamma_3 s}ds+\e M_V\int_0^te^{-\e\gamma_3 (t-s)}ds\no\\
&=&\frac{|B v_0|+M_V}{\gamma_3}(1-e^{-\e\gamma_3 t}).
\label{ches3}
\ee
Combining (\ref{ches1}) (\ref{ches2}) with (\ref{ches3}), we obtain that 
\ce
|\check{z}^{\varepsilon}_t-\tilde{\check{z}}^0_t|&\leq&\frac{e^{-\mu t}}{1-\left(\frac{L}{\gamma_1-\mu}+\frac{\e L}{\mu-\e\gamma_2}\right)}\(2|u_0-u|+2\frac{M_U}{\gamma_1}+\frac{M_V}{\gamma_2}\)\\
&&+\frac{1-\frac{\e L}{\mu-\e\gamma_2}}{1-\left(\frac{L}{\gamma_1-\mu}+\frac{\e L}{\mu-\e\gamma_2}\right)}\frac{|B v_0|+M_V}{\gamma_3}(1-e^{-\e\gamma_3 t})\\
&&+C\beta(\e).
\de
The proof is over.
\end{proof}

\br
By the estimate (\ref{malidi}), we know that when $\e$ is smaller, $\check{z}^{\varepsilon}_t$ is not nearer to $\tilde{\check{z}}^0_t$ for $t\geq0$. But this does not mean that convergence of $\check{z}^{\varepsilon}_t$ to $\tilde{\check{z}}^0_t$ fails, and it might just mean that the estimate was not good. Thus, we can not obtain the similar result to that in Theorem \ref{filcon}.
\er

\section{Conclusions}\label{con}

In the paper, we consider multiscale stochastic dynamical systems driven by L\'evy processes. First, it is proved that these systems can approximate low-dimensional systems on random invariant manifolds. Second, we establish that nonlinear filterings of multiscale stochastic dynamical systems also approximate that of reduced low-dimensional systems. Finally, we analysis the case for $\e=0$. It is unfortunate to obtain that these reduced systems does not approximate these multiscale systems.

In the future, we will investigate the possibility of doing a similar reduction and the related nonlinear filtering for the multiplicative noisy terms in the equations.

\bigskip

\textbf{Acknowledgements:}

The author would like to thank Professor Xicheng Zhang for his valuable discussions.


\begin{thebibliography}{999}

\bibitem{da} D. Applebaum: {\it L\'evy Processes and Stochastic Calculus}. Second Edition, Cambridge Univ. Press, Cambridge, 2009. 

\bibitem{la}
L. Arnold, \emph{Random Dynamical Systems}, Springer, Berlin, 1998.

\bibitem{BC}
A. Bain and D. Crisan, \emph{Fundamentals of Stochastic Filtering}, Springer, Berlin, 2009.

\bibitem{ccl} T. Caraballo, I. Chueshov, and J.A. Langa, Existence of invariant manifolds for coupled
parabolic and hyperbolic stochastic partial differential equations, {\it Nonlin.}, 18(2005)747-767.

\bibitem{ge}
G.  Evensen,
{\it Data Assimilation:The Ensemble Kalman Filter,} Springer-Verlag, Berlin, 2009.

\bibitem{Fu}  H. Fu, X. Liu and J. Duan,
Slow manifolds for multi-time-scale stochastic evolutionary systems,
\emph{Comm. Math. Sci.},
11(2013)141-162.

\bibitem{Imkeller} P. Imkeller, N. S. Namachchivaya, N. Perkowski and H. C. Yeong,
 Dimensional reduction in nonlinear filtering: a homogenization approach,
\emph{The Annals of Applied Probability},
23(2013)2290-2326.

\bibitem{iw} N. Ikeda, S. Watanabe: {\it Stochastic differential equations
and diffusion processes,} 2nd ed., North-Holland/Kodanska,
Amsterdam/Tokyo, 1989.

\bibitem{KP}
Y. Kabanov and S. Pergamenshchikov,
{\it Two-scale stochastic systems: asymptotic analysis and control,} Springer-Verlag, Berlin, 2003.

\bibitem{RG}
R. Z. Khasminskii and G.Yin,
On transition densities of singularly perturbed diffusions with fast and slow components,
\emph{SIAM J. Appl. Math.},
56(1996)1794-1819.


\bibitem{mg} L. Mitchell and G. Gottwald, Data assimilation in slow-fast systems using homogenized climate models, \emph{Journal of the Atmospheric Sciences}, 69 (2012)1359-1377.

\bibitem{Park1}  J. H. Park, N. S. Namachchivaya and H. C. Yeong,
Particle filters in a multiscale environment: Homogenized hybrid particle filter,
\emph{J. Appl. Mech.},
78(2011)1-10.

\bibitem{Park2}  J. H. Park, R. B. Sowers and N. S. Namachchivaya,
Dimensional reductionin nonlinear filtering,
\emph{ Nonlinearity},
23(2010)305-324.

\bibitem{Park3}
J. H. Park, B. Rozovskii and R. B. Sowers,
Efficient nonlinear filtering of a singularly perturbed stochastic hybrid system,
\emph{LMS Journal of Computation and Mathematics},
14(2011)254-270.

\bibitem{pa} A. Pazy, {\it Semigroups of Linear Operators and Applications to Partial Differential Equations},
Springer-Verlag, Berlin, 1983.

\bibitem{q1} H. Qiao, Stationary solutions for stochastic differential equations driven by L\'evy processes, {\it Journal of Dynamics and Differential Equations}, 29(2017)1195-1213.

\bibitem{q2} H. Qiao, Convergence of nonlinear filtering for stochastic dynamical systems with L\'evy noises, https://arxiv.org/abs/1707.07824.

\bibitem{q3} H. Qiao, Effective Filtering for Multiscale Stochastic Dynamical Systems in Hilbert Spaces, appear in \emph{J. Math. Anal. Appl}, 2020.

\bibitem{q4} H. Qiao, Y. Zhang and J. Duan, Effective filtering on a random slow manifold, {\it Nonlinearity}, 31(2018)4649-4666.

\bibitem{RBL}
B. L. Rozovskii,
\emph{Stochastic Evolution System: Linear Theory and Application to nonlinear Filtering},
Springer, New York, 1990.

\bibitem{sa} K. Sato: {\it L\'evy processes and infinitely divisible distributions,}
Cambridge University Press, 1999.

\bibitem{ZS}
Z.  Schuss,
\emph{Nonlinear Filtering and Optimal Phase Tracking},
Springer, New York, 2012.

\bibitem{Sch}
 B. Schmalfu{\ss} and  R. Schneider,
 Invariant manifolds for random dynamical systems with slow and fast variables,
\emph{J. Dyna. Diff. Equa.},
20(2008)133-164.

\bibitem{Tur}
M. Turcotte,   J. Garcia-Ojalvo,   and G. M. S\"uel,
A Genetic Timer through Noise-Induced Stabilization of an Unstable State,
\emph{Proceedings of the National Academy of Sciences of the United States of America},
41(2008)15732-7.

\bibitem{wu}
 F. Wu, T. Tian, J. B. Rawling and George Yin,
 Approximate method for stochastic chemical kinetics with two-time scales by chemical Langevin
 equations,
\emph{Journal of Chemical Physics},
144(2016)174112.

\bibitem{WR}
W. Wang and A. J. Roberts,
Slow manifold and averaging for slow-fast stochastic differential system,
\emph{J. Math. Anal. Appl},
 398(2013)822-839.

\bibitem{zhang} Y. Zhang, Z. Cheng, X. Zhang, X. Chen, J. Duan, X. Li, Data assimilation and parameter estimation for a multiscale stochastic system with $\alpha$-stable L\'evy noise. {\it Journal of Statistical Mechanics: Theory and Experiment,} (2017)113401.

\bibitem{zqd} Y. Zhang, H. Qiao and J. Duan, Effective filtering analysis for non-Gaussian dynamic systems, appear in {\it Applied Mathematics and Optimization}, 2019.
\end{thebibliography}
\end{document}